\documentclass[8pt]{article}
\usepackage{amsfonts}
\usepackage{graphicx}
\def\sinc{\mathrm{sinc}}
\input{epsf}
\usepackage{amssymb,latexsym,amsmath}
\usepackage{amsthm}
\usepackage{newlfont}

\begin{document}
\title{\textbf{Numerical solutions of Korteweg de
Vries and Korteweg de Vries-Burger's equations using computer
programming}}

\author{{Mehri
Sajjadian}$^{a,}$\thanks{m\,\underline
\,sajadian88@ms.tabrizu.ac.ir}}
\date{}
\maketitle
\begin{center}
$^a$Faculty of Mathematical Sciences, University of Tabriz, Tabriz,
Iran.
\end{center}

\begin{abstract}

\indent

In this paper, numerical and solitonic solutions of Korteweg de
Vries (KdV) and Korteweg de Vries-Burger's (KdVB) equations with
initial and boundary conditions are calculated by sinc-collocation
method. The basis of method is sinc functions. First, discretizing
time derivative of KdV and KdVB's equations using a classic finite
difference formula and space derivatives by $\theta-$ weighted
scheme between successive two time levels is applied, then Sinc
functions are used to solve these two equations. Mathematica
programming is used to solve matrix representation of these
equations. KdV equation describes behaviorof traveling waves which
is a third order non-linear partial differential equation (PDE).
Maximum absolute errors are given in Tables. The figures show
approximate solutions of these two equations. Three conservation
laws for KdV's equation are obtained.
\end{abstract}


\vskip .3cm
\section{Introduction}
First, In 1895, Korteveg and de Vries derived famous KdV equation
that describes weakly nonlinear shallow water waves and models one
directional long water wave of small amplitude , propagating in a
channel. It occurs in many field of physics such as in water waves,
plasma and fiber optics. Another example for this equation is pulse
wave propagation in blood vessels. After its discovery scientist
found solution of this equation that is called soliton. The dynamic
of solitary wave is modeled by this equation.

KdVB equation was derived by Su and Gardner [1] for a wide class of
nonlinear system in the weak non-linearity and long wavelength
approximation. The steady state solution of the KdVB equation has
been shown to model [2] weak plasma shocks propagation
perpendicularly to a magnetic field. When diffusion dominates
dispersion the steady state solutions of the KdVB equation are
monotonic shocks, and when dispersion dominates, the shocks are
oscillatory. The KdVB equation has been obtained when including
electron inertia effects in the description of weak nonlinear plasma
waves [3]. The KdVB equation has also been used in a study of wave
propagation through liquid field elastic tube [4] and for a
description of shallow water waves on viscous fluid.\\

Consider third order partial differential equations
\begin{equation}
u_t + 6 uu_x + u_{xxx} =0,\,\,\, x\in \Omega=(a,b)\subset
\mathbb{R},\,\, t>0,
\end{equation}
as KdV equation. KdV's equation has the analytical solution as
\begin{equation}
u(x,t)=0.5 sech^2 (0.5 (x-t)),
\end{equation}
and consider
\begin{equation}
u_t + \varepsilon uu_x - \nu u_{xx} + \mu u_{xxx} =0,\,\,\, x\in
\Omega=(a,b)\subset \mathbb{R},\,\, t>0,
\end{equation}
as KdVB's equation that has analytical solution given as
\begin{equation}
u(x,t)=\frac{-6 \nu^2}{25 \mu} [1+\tanh
(\frac{\nu}{10\mu}(x+\frac{6\nu^2}{25\mu} t))-\frac{1}{2} sech^2
(\frac{\nu}{10 \mu}(x+\frac{6 \nu^2}{25 \mu}t))],
\end{equation}
that space variable is defined as
\begin{equation}
x_i = a+(i-1) h,\,\, i=1,\ldots, N,\,\, h=\frac{|b-a|}{N-1}
\end{equation}
as KdVB's equation. In recent year numerous methods are used for
solving KdV and KdVB's equations. Wang Ju-Feng et al. obtained
numerical solution of the third-order nonlinear KdV equation using
the elementfree Galerkin (EFG) method which is based on the moving
least-squares approximation. A variational method is used to obtain
discrete equations, and the essential boundary conditions are
enforced by the penalty method $[5]$. The dynamics of solitary waves
is modeled by the Korteweg de Vries (KdV) equation. Jamrud Aminuddin
and Sehah, strated by discreetizing the KdV equation using the
finite difference method. The discreet form of the KdV equation is
put into a matrix form. The solution the of matrix is determined
using the Gauss-Jordan method $[6]$. Jie Shen et al. have studied
the eventual periodicity of solutions to the initial and boundary
value problem for the KdV equation on a half-line and with periodic
boundary data. They derived a representation formula for solutions
to the linearized KdV equation and rigorously establish the eventual
periodicity of these solutions $[7]$. Julio Duarte et al. employ the
Wavelet-Petrov-Galerkin method to obtain the numerical solution of
the equation Korterweg-de Vries (KdV)$[8]$. Homotopy Perturbation
Method (HPM), Variational Iteration Method (VIM) and Homotopy
Analysis Method (HAM) for the semi analytical solution of Kortweg-de
Vries (KdV) type equation are applied by Foad Saadi et al $[9]$.
Also, KdVB's equation has been solved in recent years. M. T.
Darvishi et al. have used a Numerical Solution of the KdV-Burgers'
Equation by Spectral Collocation Method and Darvishi's
Preconditionings $[10]$. In $[11]$ Riccati equation expansion method
is presented for constructing exact travelling wave solutions of
nonlinear evolution equations. The main idea of this method is to
take full advantage of the more solutions of Riccati equation to
construct exact travelling wave solutions of nonlinear evolution
equations. More new exact travelling wave solutions are obtained for
KdVB equation. On using variable transformations and proofs of
theorems, the asymptotic behaviour and the proper analytical
solution of the Korteweg-de Vries-Burgers equation have been found
in $[12]$. Anna Gao et al. studied the problem of optimal control of
the viscous KdVBs' equation. they develop a technique to utilize the
Cole-Hopf transformation to solve an optimal control problem for the
viscous KdVBs' equation $[13]$. In $[14]$ exact travelling wave and
solitary solutions for compound KdVBs equations are obtained by
using an improved sine-cosine method and the Wu elimination method.
Numerical solutions of the Korteweg-deVriesequation using the
periodic scattering transform $\mu$-representation is studied in
$[15]$. Operator Splitting Methods for Generalized Korteweg–De Vries
Equations has been discussed in $[16]$.

The paper is organized into six sections. Section 2 outlines some of
the main properties of sinc function and sinc method. In Section 3,
the discretization of KdV and KdVB equations is discussed. Section 4
outlines stability analysis and section 5 introduces errors and
conservation laws. Finally numerical results and the efficiency and
accuracy of the proposed numerical scheme is shown by considering
some numerical examples in Section 6.

\section{The Sinc function}
In this section the basis of sinc function is discussed$[17]$. The
sinc function is defined on the whole real line, $-\infty<x<\infty$,
by
\begin{equation}\label{}
    \sinc(x)=\left\{%
\begin{array}{ll}
    \frac{\sin(\pi x)}{\pi x}, & \hbox{$x\neq 0$;} \\
    1, & \hbox{$x=0$.} \\
\end{array}%
\right.
\end{equation}
For any $h>0$, the translated sinc functions with evenly spaced
nodes are given as

\begin{equation}\label{}
    S(j,h)(z)=\sinc(\frac{z-jh}{h}), \hskip .3cm j=0, \pm 1, \pm 2,
    \cdots.
\end{equation}

The sinc functions are cardinal for the interpolating points
$z_k=kh$ in the sense that
\begin{equation}\label{}
    S(j,h)(kh)=\delta_{jk}^{(0)}=\left\{%
\begin{array}{ll}
    1, & \hbox{$k=j$;} \\
    0, & \hbox{$k\neq j$.} \\
\end{array}%
\right.
\end{equation}
If $f$ is defined on the real line, then for $h>0$ the series
\begin{equation}\label{}
    C(f,h)(z)=\sum_{j=-\infty}^{\infty}f(jh)\sinc(\frac{z-jh}{h}),
\end{equation}
is called the Whittaker cardinal expansion of $f$ whenever this
series converges. They are based in the infinite strip $D_{s}$ in
the complex plane
\begin{equation}
D_{s}=\{w=u+iv: |v|<d\leq\frac{\pi}{2}\}.
\end{equation}
 Some derivatives of sinc function will be used in reduction the
 equation to matrix form so $[13]$,

\begin{equation}\label{G.3}ý
   ýI_{ji}^{(0)}=[S(j,h)(x)]\big|_{x=x_i} =ý
ý\left\{
\begin{array}{ll}ý
   1, ý& \hbox{$j=i$;} \\ý
   0, ý& \hbox{$j\neq i$,} \\ý
\end{array}%
\right.
\end{equation}ý
ý\begin{equation}\label{G.4}ý
   ýI_{ji}^{(1)}= \frac{d}{dx}[S(j,h)(x)]\big|_{x=x_i} =\frac{1}{h}ý
\left\{
\begin{array}{ll}ý
  0, ý& \hbox{$j=i$;} \\ý
  \frac{(-1)^{(k-j)}}{(k-j)}, ý& \hbox{$j\neq i$,} \\ý
\end{array}
\right.
\end{equation}ý and ý
\begin{equation}\label{G.5}ý
   ýI_{ji}^{(2)}=\frac{d^2}{d x^2}[S(j,h)(x)]\big|_{x=x_i} = \frac{1}{h^2}ý
\left\{
\begin{array}{ll}ý
   \frac{-{\pi}^2}{3}, ý& \hbox{$j=i$;} \\ý
   \frac{-2 (-1)^{(k-j)}}{(k-j)^2}, ý& \hbox{$j\neq i$.} \\ý
\end{array}
\right.
\end{equation}ý
\begin{equation}\label{G.4}ý
   ýI_{ji}^{(3)}= \frac{d}{dx}[S(j,h)(x)]\big|_{x=x_i} =\frac{1}{h^3}ý
\left\{
\begin{array}{ll}ý
  0, ý& \hbox{$j=i$;} \\ý
  \frac{(-1)^{(k-j)}}{(k-j)^3} [6 - {\pi}^2 (k-j)^2], ý& \hbox{$j\neq i$,} \\ý
\end{array}
\right.
\end{equation}
And so on, for even coefficients, where $r=1,2,\ldots$
\begin{eqnarray}ý
ý&&I_{ji}^{(2r)}=\frac{d^{2r}}{d
{x}^{2r}}[S(j,h)(x)]\big|_{x=x_i}\nonumber\\ý ý&&\nonumber\\ý ý&&=ý
ý\left\{
\begin{array}{ll}ý
   {\frac{\pi}{h}}^{2r}\frac{(-1)^r}{2r+1}, ý& \hbox{$j=i$;} \\ý
   ý\frac{(-1)^{(j-i)}}{h^{2r} (j-i)^{2r}} \sum_{l=0}^{r-1} (-1)^{l+1} \frac{2r!}{(2l+1)!} {\pi}^{2l} (iý-j)^{2l}, ý& \hbox{$j\neq i$,} \\ý
\end{array} \right.
\end{eqnarray}ý and for odd coefficients, where
 ý$r=ý1,2, ý\ldots$ý

\begin{eqnarray}ý
ý&&I_{ji}^{(2r+1)}=\frac{d^{2r}}{d {x}^{2r+1}}
[S(j,h)(x)]\big|_{x=x_i}\nonumber\\ý ý&&\nonumber\\ý ý&&=ý ý\left\{ý
\begin{array}{ll}ý
   0 ý& \hbox{$j=i$;} \\ý
   ý\frac{(-1)^{(j-i)}}{h^{2r+1}(j-i)^{2r+1}} \sum_{l=0}^{r-1} (-1)^{l} \frac{(2r+1)!}{(2l+1)!} {\pi}^{2l} (iý-j)^{2l}, ý& \hbox{$j\neq i$,} \\ý
\end{array}
\right.
\end{eqnarray}ý
\section{Survey of the method}
Consider third order partial differential equations
\begin{equation}
u_t + 6 uu_x + u_{xxx} =0,\,\,\, x\in \Omega=(a,b)\subset
\mathbb{R},\,\, t>0,
\end{equation}
as KdV equation and
\begin{equation}
u_t + \varepsilon uu_x - \nu u_{xx} + \mu u_{xxx} =0,\,\,\, x\in
\Omega=(a,b)\subset \mathbb{R},\,\, t>0,
\end{equation}
as KdVB's equation, with the initial condition
\begin{equation}
u(x,0)=f(x),\,\,\, x\in \overline{\Omega},
\end{equation}
and the boundary conditions
\begin{equation}
u(a,t)=g_a (t),\,\,\, u(b,t)=g_b (t),\,\, t\geq 0,
\end{equation}
where $\varepsilon$, $\mu$ and $\nu$ are constants.
 By discrediting time derivative of KdV's equation using a
classic finite difference formula and space derivatives by $\theta$-
weighted scheme we have ý
\begin{eqnarray}\label{G.9}
&&\frac{u^{n+1}ý - u^n}{\delta t}ý + ý\theta
(\varepsilon(uu_x)^{n+1} + \mu(u_{xxx})^{n+1})\nonumber\\
&&+ ý(1-\theta) (\varepsilon(uu_x)^ný+ \mu(u_{xxx})^n =ý0,
\end{eqnarray}
so using Taylor expansion for the term $uu_x$ and considering
$\varepsilon=6$ and $\mu=1$ we have
\begin{eqnarray}\label{G.11}
&&u^{n+1} + \delta t \theta (6[u^n u_x^{n+1}+ u_x^n u^{n+1}] +
u_{xxx}^{n+1})\nonumber\\
&&=u^n +(6(2\theta -1) (uu_x)^n + \delta t (1-\theta)u_{xxx}^n )
\end{eqnarray}ýý
for KdV equation and with same calculation for KdVB equation we have

\begin{eqnarray}
&&u^{n+1} +\delta t \theta (\varepsilon[u^n u_x^{n+1} u_x^n u^{n+1}]
- \nu u_{xx}^{n+1} + \mu u_{xxx}^{n+1})\nonumber\\
&&=u^n + \delta t \theta u^n u_x^n - \delta t (1-\theta)
[\varepsilon uu_x^n - \nu u_{xx}^n + \mu u_{xxx}^n]
\end{eqnarray}
Now we use approximate solution as
\begin{equation}
u(x,t^n)=u^n(x) \simeq \sum_{j=1}^N u^n_j S_j (x).
\end{equation}
where
\begin{equation}
S_j(x) = Sin(\frac{x-(j-1)h-a}{h}).
\end{equation}
By substituting above approximate solution in Eq. (22) a matrix
representation is obtained for KdV equation as
\begin{equation}
Mu^{n+1}=R
\end{equation}
where
\begin{eqnarray*}
&&A_d=[I_{ij}^{(0)}: i=ý2,..., N-1, ýj=1,..., ýN , ý\,\, 0\,\, elsewhere]_{N\times N},\\
&&A_b=[I_{ij}^{(0)}: ýi=ý1,N,j=ý1,\ldots, ýNý, ý\,\, 0\,\, elsewhere]_{N\times ýN},\\
&&B=-[I_{ij}^{(1)}: ýi=ý2,\ldots, ýNý-1, j=ý1,\ldots, ýNý, ý\,\, 0\,\,elsewhere]_{N\times ýN},\\
&&C=[I_{ij}^{(2)}: ýi=ý2,\ldots,N-1,ýj=ý1,\ldots, ýNý, ý\,\, 0\,\,elsewhere]_{N\times ýN},\\
&&G=-[I_{ij}^{(3)}: ýi=ý2,\ldots,N-1,j=ý1,\ldots, ýNý, ý\,\, 0\,\,elsewhere]_{N\times ýN},\\
&&u_x^n= B u^ný, ý\;\;\; D=p(u^n)^{p-1} * u_x^n * A_d,\;\;\; E=(u^n)*B,\\
&&F^{n+1}=[g_a (t^{ný+1}), 0, ý\ldotsý, 0, ýg_b(t^{ný+1})]^T,
\end{eqnarray*}
so with these definition for KdV equation we have
\begin{eqnarray*}
&&M=[A_dý + ýA_bý+ý\theta \delta t(6(E+D)+ G)],\\
&&R=[A_dý - \delta t \{ý6(2 \theta -1) Eý -
ý(1-\theta) G\}] u^ný +F^{n+1},ý\\
\end{eqnarray*}
and with same substitution for KdVB's equation
\begin{eqnarray*}
&&M=[A_dý + ýA_bý+ý\theta \delta (\varepsilon (E+D)+ \mu G - \nu C)],\\
&&R=[A_dý + \delta t \{ý\varepsilon(2 \theta -1) Eý -
ý(1-\theta) (\mu G -\nu C)\}] u^ný +F^{n+1},ý\\
\end{eqnarray*}
\section{Stability Analysis}
In this section stability analysis of approximate solution for
linearized equation is discussed[18]. The error at nth time level is
$$e^n = u^n_{exact} - u^n_{approximate}.$$
\subsection{KdV's equation}
By considering the obtained matrix we have
\begin{equation}
[H+\delta t \theta K] e^{n+1} = [H -\delta t (1- \theta )K] e^n,
\end{equation}
where $H=[A_d + A_b]A^{-1}$ and $K=[6 E + G]A^{-1}$ so,
$$e^{n+1} = P e^n$$
where $P=[H+\delta t \theta K]^{-1}[H -\delta t (1- \theta )K]$.
This method is stable if $\parallel P \parallel_2 \leq 1$ or $\rho
(P) \leq 1$ which is spectral radius of the matrix $P$. The
stability is assured if all the eigenvalues of the matrix $[H+\delta
t \theta K]^{-1}[H -\delta t (1- \theta )K]$ satisfy the following
condition
\begin{equation}
\bigg|\frac{\lambda_H -\delta t (1- \theta ) \lambda_K}{\lambda_H +
\delta t \theta \lambda_K}\bigg|\leq 1
\end{equation}
where ý${\lambda}_H$ and ${\lambda}_K$ýare eigenvalues of the
matrices $H$ and $K$ respectively. When $\theta=0.5$, the inequality
(28) becomes
\begin{equation}\label{G.19}
\bigg|\frac{{\lambda}_Hý - 0.5 ý\delta t {\lambda}_K}{{\lambda}_Hý +
0.5 ý\delta t {\lambda}_K}\bigg|\leqý 1.
\end{equation}ý
In the case of complex eigenvalues ý${\lambda}_H = a_hý + ýib_h$ýand
${\lambda}_K = a_ký + ýib_k$, where $a_h$, $a_k$, $b_h$ and $b_k$
are any real numbers, the inequality (29) takes the following form,
\begin{equation}\label{G.20}ý
ý\bigg|\frac{({a}_hý - 0.5 ý\delta t {a}_k)ý + ýi(b_hý-0.5 ý\delta t
b_k)}{({a}_hý - 0.5 ý\delta t {a}_k)ý + ýi(b_hý-0.5 ý\delta t
b_k)}\bigg|\leqý 1.
\end{equation}ý
The inequality (30) is satisfied if $a_h a_k + b_h b_k \geq 0$. For
real eigenvalues, the inequality (29) holds true if either
($\lambda_h \geq 0$ and $\lambda_k\geq 0$) or ($\lambda_h \leq 0$
and $\lambda_k\leq 0$). This shows that the scheme is
unconditionally stable if $a_h a_k + b_h b_k \geq 0$, for complex
eigenvalues and if either ($\lambda_h \geq 0$ and $\lambda_k\geq 0$)
or ($\lambda_h \leq 0$ and $\lambda_k\leq 0$), for real eigenvalues.
When $\theta = 0$, the inequality (29) becomes
\begin{equation*}
\bigg|1-\frac{\delta t \lambda_k}{\lambda_h}\bigg|\leq 1,
\end{equation*}
i.e.,
$$\delta t \leq \frac{2\lambda_H}{\lambda_K}\,\,\, and\,\,\, \frac{\lambda_H}{\lambda_K}\geq 0$$
Thus for $\theta = 0$, the scheme is conditionally stable. The
stability of scheme for the other values of è can be investigate in
a similar manner. The stability of the scheme and conditioning of
the component matrices H, K of the matrix P depend on the weight
parameter è and the minimum distance between any two collocation
points $h$ in the domain set $[a, b]$.
\subsection{KdVB's equation}
By considering the obtained matrix we have
\begin{equation}
[H+\delta t \theta K] e^{n+1} = [H -\delta t (1- \theta )K] e^n,
\end{equation}
where $H=[A_d + A_b]A^{-1}$ and $K=[\varepsilon E + \mu G -\nu
C]A^{-1}$ so,
$$e^{n+1} = P e^n$$
where $P=[H+\delta t \theta K]^{-1}[H -\delta t (1- \theta )K]$.
This method is stable if $\parallel P \parallel_2 \leq 1$ or $\rho
(P) \leq 1$ which is spectral radius of the matrix $P$. The
stability is assured if all the eigenvalues of the matrix $[H+\delta
t \theta K]^{-1}[H -\delta t (1- \theta )K]$ satisfy the following
condition
\begin{equation}
\bigg|\frac{\lambda_H -\delta t (1- \theta ) \lambda_K}{\lambda_H +
\delta t \theta \lambda_K}\bigg|\leq 1
\end{equation}
where ý${\lambda}_H$ and ${\lambda}_K$ýare eigenvalues of the
matrices $H$ and $K$ respectively. When $\theta=0.5$, the inequality
(32) becomes
\begin{equation}\label{G.19}
\bigg|\frac{{\lambda}_Hý - 0.5 ý\delta t {\lambda}_K}{{\lambda}_Hý +
0.5 ý\delta t {\lambda}_K}\bigg|\leqý 1.
\end{equation}ý
In the case of complex eigenvalues ý${\lambda}_H = a_hý + ýib_h$ýand
${\lambda}_K = a_ký + ýib_k$, where $a_h$, $a_k$, $b_h$ and $b_k$
are any real numbers, the inequality (33) takes the following form,
\begin{equation}\label{G.20}ý
ý\bigg|\frac{({a}_hý - 0.5 ý\delta t {a}_k)ý + ýi(b_hý-0.5 ý\delta t
b_k)}{({a}_hý - 0.5 ý\delta t {a}_k)ý + ýi(b_hý-0.5 ý\delta t
b_k)}\bigg|\leqý 1.
\end{equation}ý
The inequality (34) is satisfied if $a_h a_k + b_h b_k \geq 0$. For
real eigenvalues, the inequality (33) holds true if either
($\lambda_h \geq 0$ and $\lambda_k\geq 0$) or ($\lambda_h \leq 0$
and $\lambda_k\leq 0$). This shows that the scheme is
unconditionally stable if $a_h a_k + b_h b_k \geq 0$, for complex
eigenvalues and if either ($\lambda_h \geq 0$ and $\lambda_k\geq 0$)
or ($\lambda_h \leq 0$ and $\lambda_k\leq 0$), for real eigenvalues.
When $\theta = 0$, the inequality (33) becomes
\begin{equation*}
\bigg|1-\frac{\delta t \lambda_k}{\lambda_h}\bigg|\leq 1,
\end{equation*}
i.e.,
$$\delta t \leq \frac{2\lambda_H}{\lambda_K}\,\,\, and\,\,\, \frac{\lambda_H}{\lambda_K}\geq 0$$
Thus for $\theta = 0$, the scheme is conditionally stable. The
stability of scheme for the other values of è can be investigate in
a similar manner. The stability of the scheme and conditioning of
the component matrices H, K of the matrix P depend on the weight
parameter è and the minimum distance between any two collocation
points $h$ in the domain set $[a, b]$.

\section{Errors and Conservation Laws}
In this section, two error norms is defined that will be used for
showing the accuracy of the method as follows
\begin{eqnarray*}
&&L_2 = \parallel u - \tilde{u} \parallel_2 = \sqrt{h \sum_{j=1}^{N}|u - \tilde{u_{j}}|^2},\\
&&L_{\infty} = \parallel u - \tilde{u} \parallel_{\infty} = {max}
|u_j -\tilde{u_j}|,\,\, {1\leq j \leq N}
\end{eqnarray*}
where $u,\,\,\tilde{u}$ are exact and approximate solution
respectively. \\
KdV equation has three conservation laws as follows
\begin{eqnarray*}
&&I_1 = \int_a^b u(x,t) dx,\\
&&I_2=\int_a^b u(x,t)^2 dx,\\
&&I_3=\int_a^b [u(x,t)^2 -\frac{1}{3} u(x,t)^3] dx,
\end{eqnarray*}
Where $I_1,\,\, I_2, \,\, I_3$ represents mass, momentum and energy
that show conservation properties of collocation method by,
\begin{eqnarray*}
&&I_1 \simeq h \sum_{j=1}^N u_j^n,\\
&&I_2 \simeq h \sum_{j=1}^N (u_j^n)^2,\\
&&I_3 \simeq h \sum_{j=1}^N [(u_j^n)^2 -\frac{1}{3} (u_j^n)^3],
\end{eqnarray*}
Sinc-collcation scheme satisfies the properties $I_1$ and $I_2$ in
the
conservative case.\\
\section{Numerical Solution}
In this section $L_2$ and $L_{\infty}$ are obtained and shown in
Tables and approximate solution of KdV and KdVB's equations is shown
in Figures.

\subsection{KdV equaiton}
By considering $\mu=1,\, \varepsilon =6,\, \theta=0.5$, in Table 1,
two kinds of error is calculated for $n=100,\,a=-15,\,b=15,\, \delta
t=0.1,\, T=0.1,\ldots, 0.9$. Three invariants of conservative laws
is obtained. Table 2 indicates errors for $\delta t=0.01$ and Table
3, shows error for $delta t =0.001$. Table 4, errors are given for
different time levels $T=0.001,\ldots,0.009$. By Tables, Errors
reduce when time step decreases. Figure 1, indicates solitonic
solutions of KdV equation in different time level of $T=1,\ldots,9$
for $n=100,\, a=-10,\, b=20,\, \delta t =0.01$. Figure 2, indicates
solitonic solutions of KdV equation in different time level of
$T=1,\ldots,9$ for $n=100,\, a=-15,\, b=15,\, \delta t=0.1$.
\begin{center}
  \begin{tabular}{||c|c|c|c|c|c||}
    \hline \hline
    $Time$ & $L_{\infty}$ & $L_2$ & $I_1$ & $I_2$ & $I_3$\\
    \hline
    0.1 &4.55798$\times {10}^{-5}$ & 7.31340$\times {10}^{-5}$ & 2  & 0.666667  & -0.0888889\\
    0.2& 9.28126$\times {10}^{-5}$ & 1.34140$\times {10}^{-4}$  & 2  & 0.666667  & -0.0888853\\
    0.3 &1.23469$\times {10}^{-4}$ & 1.81748$\times {10}^{-4}$  & 2  & 0.666667  & -0.0888766\\
    0.4 &1.34961$\times {10}^{-4}$ & 2.20735$\times {10}^{-4}$  & 2  & 0.666667  & -0.0888665\\
    0.5 &1.39977$\times {10}^{-4}$ & 2.54695$\times {10}^{-4}$  & 2  & 0.666667  & -0.0888567\\
    0.6 &1.58066$\times {10}^{-4}$ & 2.85701$\times {10}^{-4}$  & 2  & 0.666667  & -0.0888476\\
    0.7 &1.77771$\times {10}^{-4}$ & 3.15070$\times {10}^{-4}$  & 2  & 0.666667  & -0.0888394\\
    0.8 &1.98187$\times {10}^{-4}$ & 3.43619$\times {10}^{-4}$  & 2.00001  & 0.666667  & -0.0888321\\
    0.9 &2.20096$\times {10}^{-4}$ & 3.71849$\times {10}^{-4}$  & 2.00001  & 0.666667  & -0.0888257\\
    \hline \hline
  \end{tabular}\\
  \vskip .3cm
\end{center}
\begin{center}
Table1: Errors and invariants for $n=100,\,a=-15,\,b=15,\, \delta
t=0.1,\, T=0.1,\ldots, 0.9$.
\end{center}
\begin{center}
  \begin{tabular}{||c|c|c|c|c|c||}
    \hline \hline
    $Time$ & $L_{\infty}$ & $L_2$ & $I_1$ & $I_2$ & $I_3$\\
    \hline
    0.1 &1.74342$\times {10}^{-6}$ & 1.91251$\times {10}^{-6}$ & 2  & 0.666667  & -0.0888889\\
    0.2& 2.31194$\times {10}^{-6}$ & 2.87693$\times {10}^{-6}$  & 2  & 0.666667  & -0.0888888\\
    0.3 &2.73488$\times {10}^{-6}$ & 3.64966$\times {10}^{-6}$  & 2  & 0.666667  & -0.0888887\\
    0.4 &3.04988$\times {10}^{-6}$ & 4.32025$\times {10}^{-6}$  & 2  & 0.666667  & -0.0888886\\
    0.5 &3.32879$\times {10}^{-6}$ & 4.95178$\times {10}^{-6}$  & 2  & 0.666667  & -0.0888885\\
    0.6 &3.49763$\times {10}^{-6}$ & 5.53733$\times {10}^{-6}$  & 2  & 0.666667  & -0.0888884\\
    0.7 &3.93719$\times {10}^{-6}$ & 6.16498$\times {10}^{-6}$  & 2  & 0.666667  & -0.0888883\\
    0.8 &1.98187$\times {10}^{-6}$ & 6.75605$\times {10}^{-6}$  & 2.00001  & 0.666667  & -0.0888883\\
    0.9 &4.69034$\times {10}^{-6}$ & 7.26574$\times {10}^{-6}$  & 2.00001  & 0.666667  & -0.0888882\\
    \hline \hline
  \end{tabular}\\
  \vskip .3cm
\end{center}
\begin{center}
Table2: Errors and invariants for $n=100,\,a=-15,\,b=15,\, \delta t
=0.01,\, T=0.1,\ldots, 0.9$.
\end{center}

\begin{center}
  \begin{tabular}{||c|c|c|c|c|c||}
    \hline \hline
    $Time$ & $L_{\infty}$ & $L_2$ & $I_1$ & $I_2$ & $I_3$\\
    \hline
    0.1 &1.74199$\times {10}^{-6}$ & 1.75980$\times {10}^{-6}$  & 2  & 0.666667  & -0.0888889\\
    0.2& 2.30437$\times {10}^{-6}$ & 2.54006$\times {10}^{-6}$  & 2  & 0.666667  & -0.0888889\\
    0.3 &2.71876$\times {10}^{-6}$ & 3.16315$\times {10}^{-6}$  & 2  & 0.666667  & -0.0888889\\
    0.4 &3.07014$\times {10}^{-6}$ & 3.70810$\times {10}^{-6}$  & 2  & 0.666667  & -0.0888889\\
    0.5 &3.27326$\times {10}^{-6}$ & 4.23733$\times {10}^{-6}$  & 2  & 0.666667  & -0.0888889\\
    0.6 &3.69103$\times {10}^{-6}$ & 4.73099$\times {10}^{-6}$  & 2  & 0.666667  & -0.0888889\\
    0.7 &4.03713$\times {10}^{-6}$ & 5.23450$\times {10}^{-6}$  & 2  & 0.666667  & -0.0888889\\
    0.8 &4.20408$\times {10}^{-6}$ & 5.77575$\times {10}^{-6}$  & 2  & 0.666667  & -0.0888889\\
    0.9 &4.59416$\times {10}^{-6}$ & 6.26060$\times {10}^{-6}$  & 2  & 0.666667  & -0.0888889\\
    \hline \hline
  \end{tabular}\\
  \vskip .3cm
\end{center}
\begin{center}
Table3: Errors and invariants for $n=100,\,a=-15,\,b=15,\, \delta t
=0.001,\, T=0.1,\ldots, 0.9$.
\end{center}

\begin{center}
  \begin{tabular}{||c|c|c|c|c|c||}
    \hline \hline
    $Time$ & $L_{\infty}$ & $L_2$ & $I_1$ & $I_2$ & $I_3$\\
    \hline
    0.01 &7.90227$\times {10}^{-7}$ & 4.72378$\times {10}^{-7}$  & 2  & 0.666667  & -0.0888889\\
    0.02& 1.03034$\times {10}^{-6}$ & 7.25023$\times {10}^{-7}$  & 2  & 0.666667  & -0.0888889\\
    0.03 &1.14705$\times {10}^{-6}$ & 9.14997$\times {10}^{-7}$  & 2  & 0.666667  & -0.0888889\\
    0.04 &1.21740$\times {10}^{-6}$ & 1.07375$\times {10}^{-6}$  & 2  & 0.666667  & -0.0888889\\
    0.05 &1.28565$\times {10}^{-6}$ & 1.21327$\times {10}^{-6}$  & 2  & 0.666667  & -0.0888889\\
    0.06 &1.41361$\times {10}^{-6}$ & 1.33912$\times {10}^{-6}$  & 2  & 0.666667  & -0.0888889\\
    0.07 &1.51656$\times {10}^{-6}$ & 1.45490$\times {10}^{-6}$  & 2  & 0.666667  & -0.0888889\\
    0.08 &1.60373$\times {10}^{-6}$ & 1.56261$\times {10}^{-6}$  & 2  & 0.666667  & -0.0888889\\
    0.09 &1.67755$\times {10}^{-6}$ & 1.66392$\times {10}^{-6}$  & 2  & 0.666667  & -0.0888889\\
    \hline \hline
  \end{tabular}\\
  \vskip .3cm
\end{center}
\begin{center}
Table4:Errors and invariants for $n=100,\,a=-15,\,b=15,\, \delta t
=0.001,\, T=0.1,\ldots, 0.9 and T=0.001,\ldots, 0.009$.
\end{center}
\newpage
{\setlength{\unitlength}{1mm}
\begin{picture}(10,10)
\includegraphics{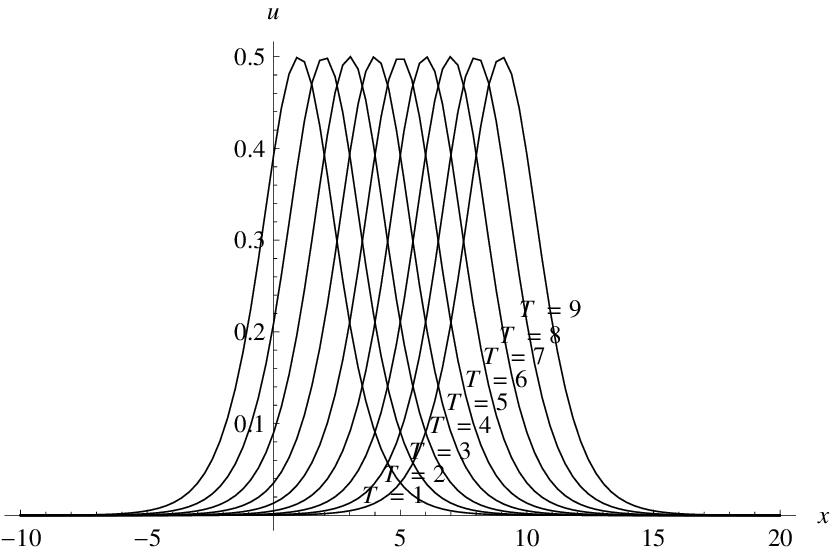}
\end{picture}\\
\vskip 7cm

Fig1: Solitonic solutions of KdV equation in different time level of
$T=1,\ldots,9$ for $n=100,\, a=-10,\, b=20,\, \delta t =0.01$

{\setlength{\unitlength}{1mm}
\begin{picture}(10,10)
\includegraphics{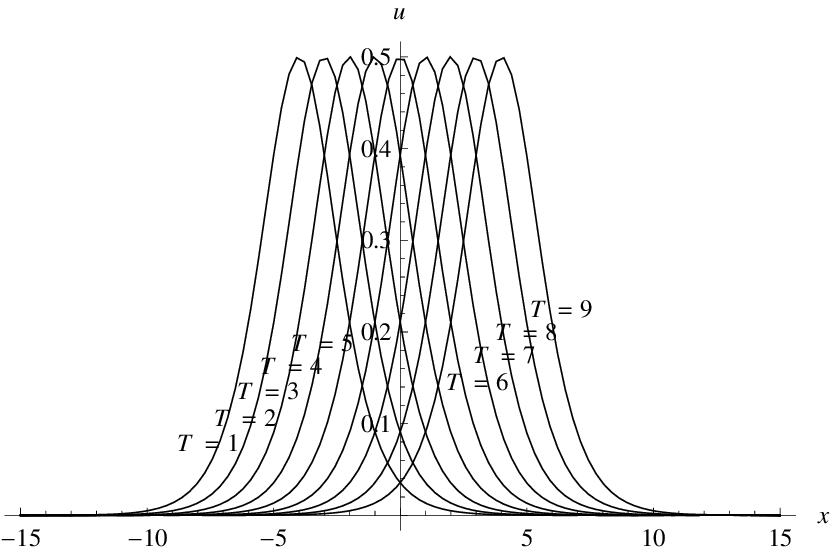}
\end{picture}\\
\vskip 7cm

Fig2: Solitonic solutions of KdV equation in different time level of
$T=1,\ldots,9$ for $n=100,\, a=-15,\, b=15,\, \delta t =0.1$

\subsection{KdVB equation}
By considering $\mu=0.1,\,\varepsilon =2, \nu=0.005,\,\theta=0.5$,
in Table 5, two kinds of error is calculated for
$n=100,\,a=-100,\,b=100,\, \delta t=0.02,\, T=1,\ldots, 9$.

Table 6 indicates errors for $\mu=0.001,\,\varepsilon =1,
\nu=0.001,\,\theta=0.5$ and $n=100,\,a=-40,\,b=100,\, \delta
t=0.05,\, T=1,\ldots, 9$.

In Table 7, error is obtained for $n=16,\,a=0,\,b=100,\, \delta
t=0.00001,\, T=0.0001,\ldots, 0.0009$ and $\mu=0.001,\,\varepsilon
=1, \nu=0.001,\,\theta=0.5$.

Table 8 shows errors for $n=16,\,a=-40,\,b=40,\, \delta t=0.02,\,
T=1,\ldots, 9$ and $\mu=0.1,\,\varepsilon =2,
\nu=0.005,\,\theta=0.5$.

In Table 9, errors are calculated for $n=16,\,a=8,\,b=99,\, \delta
t=0.05,\, T=1,\ldots, 9$ and $\mu=0.001,\,\varepsilon =1,
\nu=0.001,\,\theta=0.5$. Errors reduces by decreasing of time steps.

Figure 3, indicates solutions of KdVB equation for $n=100,\,
a=-100,\, b=100,\, \delta t=0.02,\, T=1$ and $\mu=0.1,\,\varepsilon
=2, \nu=0.005,\,\theta=0.5$.

Figure 4, indicates solutions of KdVB equation for $n=100,\,
a=-40,\, b=100,\, \delta t=0.05,\, T=1$ and $\mu=0.1,\,\varepsilon
=1, \nu=0.1,\,\theta=0.5$.

Figure 5 shows solutions of KdVB equation for $n=100,\, a=-40,\,
b=100,\, \delta t=0.05,\, T=1$ and $\mu=0.01,\,\varepsilon =1,
\nu=0.01,\,\theta=0.5$.

In Figure 6, numerical solution of KdVB equation is shown for
$n=100,\, a=-40,\, b=100,\, \delta t=0.05,\, T=1$ and
$\mu=0.001,\,\varepsilon =1, \nu=0.001,\,\theta=0.5$

\begin{center}
  \begin{tabular}{||c|c|c||}
    \hline \hline
    $Time$ & $L_{\infty}$ & $L_2$ \\
    \hline
   1 &2.06127$\times {10}^{-6}$ & 4.62479$\times {10}^{-6}$  \\
   2& 1.03034$\times {10}^{-6}$ & 7.25023$\times {10}^{-7}$  \\
   3 &1.14705$\times {10}^{-6}$ & 9.14997$\times {10}^{-7}$  \\
   4 &1.21740$\times {10}^{-6}$ & 1.07375$\times {10}^{-6}$  \\
   5 &1.28565$\times {10}^{-6}$ & 1.21327$\times {10}^{-6}$  \\
   6 &1.41361$\times {10}^{-6}$ & 1.33912$\times {10}^{-6}$  \\
   7 &1.51656$\times {10}^{-6}$ & 1.45490$\times {10}^{-6}$  \\
   8 &1.60373$\times {10}^{-6}$ & 1.56261$\times {10}^{-6}$  \\
   9 &1.67755$\times {10}^{-6}$ & 1.66392$\times {10}^{-6}$  \\
    \hline \hline
\end{tabular}\\
Table5: Errors for $\mu=0.1,\,\varepsilon =2,
\nu=0.005,\,\theta=0.5$ and $n=100,\,a=-100,\,b=100,\, \delta
t=0.02,\, T=1,\ldots, 9$.
\end{center}

\begin{center}
  \begin{tabular}{||c|c|c||}
    \hline \hline
    $Time$ & $L_{\infty}$ & $L_2$ \\
    \hline
   1 &5.49372$\times {10}^{-7}$ & 1.00098$\times {10}^{-6}$  \\
   2& 1.09951$\times {10}^{-6}$ & 1.99828$\times {10}^{-6}$  \\
   3 &1.65039$\times {10}^{-6}$ & 2.99190$\times {10}^{-6}$  \\
   4 &2.20200$\times {10}^{-6}$ & 3.98187$\times {10}^{-6}$  \\
   5 &2.75434$\times {10}^{-6}$ & 4.96819$\times {10}^{-6}$  \\
   6 &3.30738$\times {10}^{-6}$ & 5.95089$\times {10}^{-6}$  \\
   7 &3.86111$\times {10}^{-6}$ & 6.92997$\times {10}^{-6}$  \\
   8 &4.41551$\times {10}^{-6}$ & 7.90547$\times {10}^{-6}$  \\
   9 &4.97059$\times {10}^{-6}$ & 8.87738$\times {10}^{-6}$  \\
    \hline \hline
\end{tabular}\\
\vskip .3cm

Table6: Errors for $\mu=0.001,\,\varepsilon =1,
\nu=0.001,\,\theta=0.5$ and $n=100,\,a=-40,\,b=100,\, \delta
t=0.05,\, T=1,\ldots, 9$.
\end{center}

\begin{center}
  \begin{tabular}{||c|c|c||}
    \hline \hline
    $Time$ & $L_{\infty}$ & $L_2$ \\
    \hline
   1 &1.84693$\times {10}^{-12}$ & 6.66561$\times {10}^{-12}$  \\
   2& 3.69387$\times {10}^{-12}$ & 1.33287$\times {10}^{-11}$  \\
   3 &5.54080$\times {10}^{-12}$ & 1.99924$\times {10}^{-11}$  \\
   4 &7.38774$\times {10}^{-12}$ & 2.66562$\times {10}^{-11}$  \\
   5 &9.23467$\times {10}^{-12}$ & 3.33201$\times {10}^{-11}$  \\
   6 &1.10816$\times {10}^{-11}$ & 3.99840$\times {10}^{-11}$  \\
   7 &1.29285$\times {10}^{-11}$ & 4.66479$\times {10}^{-11}$  \\
   8 &1.47755$\times {10}^{-11}$ & 5.33118$\times {10}^{-11}$  \\
   9 &1.66224$\times {10}^{-11}$ & 5.99758$\times {10}^{-11}$  \\
    \hline \hline
\end{tabular}\\
\vskip .3cm

Table7: Errors for $\mu=0.001,\,\varepsilon =1,
\nu=0.001,\,\theta=0.5$ and $n=16,\,a=0,\,b=100,\, \delta
t=0.00001,\, T=0.0001,\ldots, 0.0009$.
\end{center}

\begin{center}
  \begin{tabular}{||c|c|c||}
    \hline \hline
    $Time$ & $L_{\infty}$ & $L_2$ \\
    \hline
   1 &7.82718$\times {10}^{-8}$ & 3.31435$\times {10}^{-7}$  \\
   2& 1.57051$\times {10}^{-7}$ & 6.62416$\times {10}^{-7}$  \\
   3 &2.36335$\times {10}^{-7}$ & 9.92942$\times {10}^{-7}$  \\
   4 &3.16123$\times {10}^{-7}$ & 1.32301$\times {10}^{-6}$  \\
   5 &3.96412$\times {10}^{-7}$ & 1.65261$\times {10}^{-6}$  \\
   6 &4.77201$\times {10}^{-7}$ & 1.98176$\times {10}^{-6}$  \\
   7 &5.58486$\times {10}^{-7}$ & 2.31043$\times {10}^{-6}$  \\
   8 &6.40266$\times {10}^{-7}$ & 2.63863$\times {10}^{-6}$  \\
   9 &7.22539$\times {10}^{-7}$ & 2.96637$\times {10}^{-6}$  \\
    \hline \hline
\end{tabular}\\
\vskip .3cm

Table8: Errors for $\mu=0.1,\,\varepsilon =2,
\nu=0.005,\,\theta=0.5$ and $n=16,\,a=-40,\,b=40,\, \delta t=0.02,\,
T=1,\ldots, 9$.
\end{center}

\begin{center}
  \begin{tabular}{||c|c|c||}
    \hline \hline
    $Time$ & $L_{\infty}$ & $L_2$ \\
    \hline
   1 &2.26009$\times {10}^{-8}$ & 8.83534$\times {10}^{-8}$  \\
   2& 4.52028$\times {10}^{-8}$ & 1.76687$\times {10}^{-7}$  \\
   3 &6.78055$\times {10}^{-8}$ & 2.65005$\times {10}^{-7}$  \\
   4 &9.04092$\times {10}^{-8}$ & 3.53306$\times {10}^{-7}$  \\
   5 &1.13014$\times {10}^{-7}$ & 4.41592$\times {10}^{-7}$  \\
   6 &1.35619$\times {10}^{-7}$ & 5.29861$\times {10}^{-7}$  \\
   7 &1.58226$\times {10}^{-7}$ & 6.18113$\times {10}^{-7}$  \\
   8 &1.80833$\times {10}^{-7}$ & 7.06350$\times {10}^{-7}$  \\
   9 &2.03441$\times {10}^{-7}$ & 7.94570$\times {10}^{-7}$  \\
    \hline \hline
\end{tabular}\\
\vskip .3cm

Table9: Errors for $\mu=0.001,\,\varepsilon =1,
\nu=0.001,\,\theta=0.5$ and $n=16,\,a=8,\,b=99,\, \delta t=0.05,\,
T=1,\ldots, 9$.
\end{center}

\newpage
{\setlength{\unitlength}{1mm}
\begin{picture}(10,10)
\includegraphics{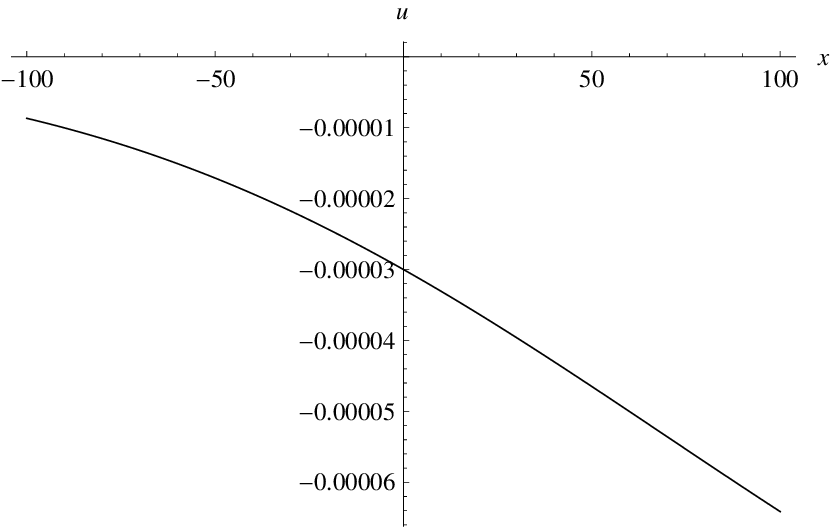}
\end{picture}\\
\vskip 7cm

Fig3: Solutions of KdVB equation for $n=100,\, a=-100,\, b=100,\,
\delta t=0.02,\, T=1$ and $\mu=0.1,\,\varepsilon =2,
\nu=0.005,\,\theta=0.5$

{\setlength{\unitlength}{1mm}
\begin{picture}(10,10)
\includegraphics{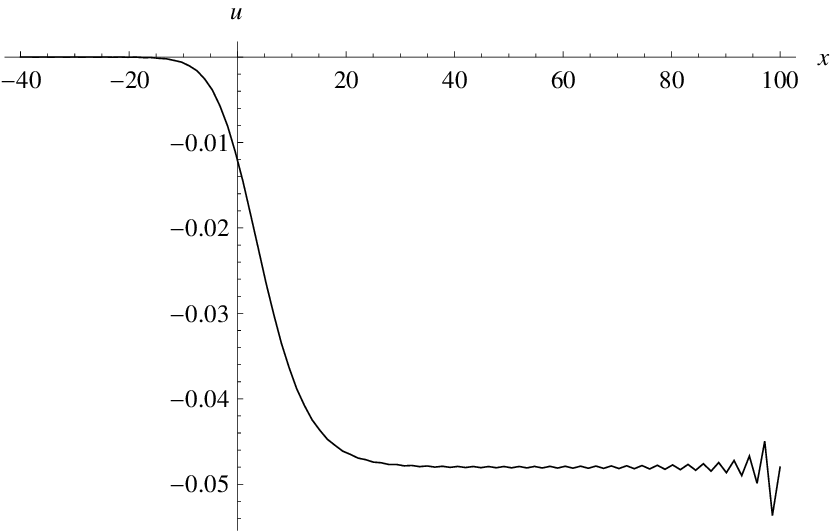}
\end{picture}\\
\vskip 7cm

Fig4: solutions of KdVB equation for $n=100,\, a=-40,\, b=100,\,
\delta t=0.05,\, T=1$ and $\mu=0.1,\,\varepsilon =1,
\nu=0.1,\,\theta=0.5$.

\newpage
{\setlength{\unitlength}{1mm}
\begin{picture}(10,10)
\includegraphics{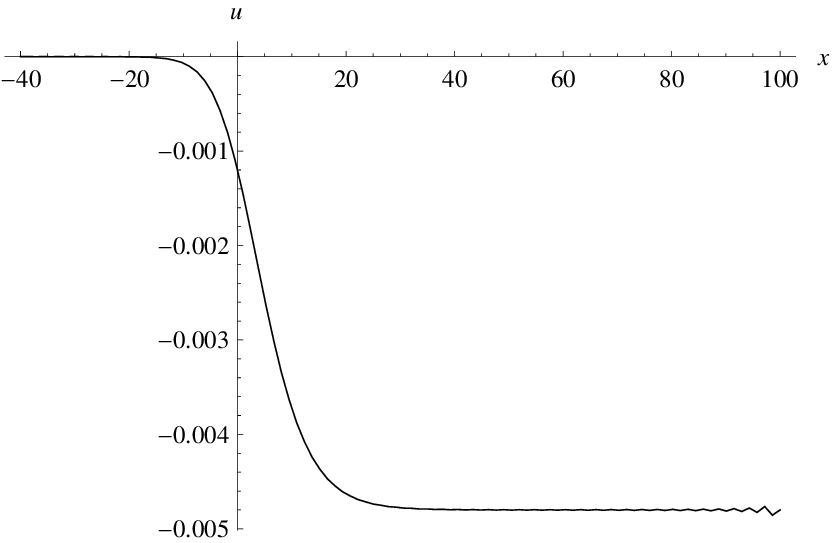}
\end{picture}\\
\vskip 7cm

Fig5:solutions of KdVB equation for $n=100,\, a=-40,\, b=100,\,
\delta t=0.05,\, T=1$ and $\mu=0.01,\,\varepsilon =1,
\nu=0.01,\,\theta=0.5$.

{\setlength{\unitlength}{1mm}
\begin{picture}(10,10)
\includegraphics{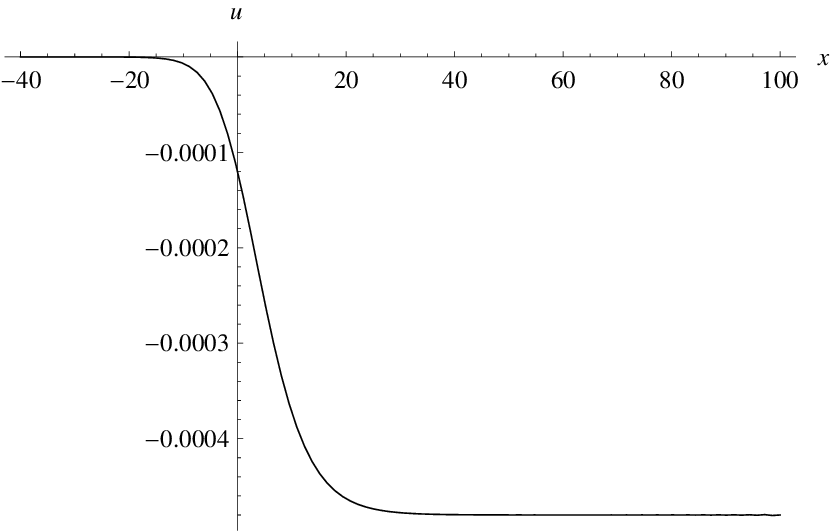}
\end{picture}\\
\vskip 7cm

Fig6: Solution of KdVB equation for $n=100,\, a=-40,\, b=100,\,
\delta t=0.05,\, T=1$ and $\mu=0.001,\,\varepsilon =1,
\nu=0.001,\,\theta=0.5$

\section*{Conclusion}
The collocation method using sinc basis is applied for solving KdV
and KdVB's equations. Three invariant of conservation laws are
calculated for KdV equation. The method is computationally
attractive and results are shown through tables and figures.


\section*{References}
\begin{description}

   \item[[1]]. C.H. Su and C.s. Gardner," Drivation of the Korteweg-de Vries
   and Burgers' Equation ", J. Math. Phy., 10, 536-539, (1969).

   \item[[2]]. H. Gard and P. N. Hu, " Unified shock profile in plasma", Phys.
   Fluids, 10, 2596- 2602, (1967).

   \item[[3]]. H. Grad and P. N. Hu, "
   Collisional theory of shock and nonlinear waves in plasma", J. phys.
   Fluids, 15, 845-864, (1972).

   \item[[4]]. R. S. johnson, " A nonlinear
   equation incorporating damping and dispersion", J. Phys. Fluids
   Mech. 42, 94-60, (1970).

   \item[[5]]. Wang Ju-Feng, Sun Feng-Xin, and Cheng Rong-Jun, " Element-free Galerkin method for a kind of KdV equation" , Chin.
   Phys. B, 19, 6, (2010).

   \item[[6]]. Jamrud Aminuddin, and Sehah, "Numerical Solution of The Korteweg de Vries Equation", Int. J. of Basic and Appl. Sciences IJBAS-IJENS,
   11, 2, (2011).

   \item[[7]]. Jie Shen, Jiahong Wu, Juan-Ming Yuan, "Eventual periodicity for the KdV equation on a
   half-line", Physica D, 227, 105–119 (2007).

   \item[[8]]. Jairo Villegas G. and Jorge Casta$\tilde{n}$o B. "Wavelet-Petrov-Galerkin Method for the
   Numerical Solution of the KdV Equation", Appl. Math. Sci, 6, 69, 3411 - 3423, (2012).

   \item[[9]]. Foad Saadi, M. Jalali Azizpour and S.A. Zahedi,
   "Analytical Solutions of Kortweg-de Vries
   (KdV) Equation", World Academy of Science, Engineering and
   Technology, 69, (2010).

   \item[[10]]. M. T. Darvishia, F. Khanib and S. Kheybaria, "Numerical Solution of the KdV-Burgers' Equation
   by Spectral Collocation Method and Darvishi's Preconditionings",
   Int. J. Contemp. Math. Sciences, 2, 22, 1085 - 1095, (2007).

   \item[[11]]. Yang Feng, Hong-qing Zhang, "Riccati Equation Expansion Method for Solving KdV-Burgers Equation", Int. J. of Nonlinear
   Sci., 7, 1, 84-89, (2009).

   \item[[12]]. Shu Jian-Jun, "The proper analytical solution of the
   Korteweg-de Vries-Burgers equation", J. Phys. A: Math. Gen., 20, L49-L56, (1987).

   \item[[13]]. Anna Gao, Chunyu Shen, Xinghua Fan, \textit{Optimal Control of the Viscous KdV-Burgers Equation Using an Equivalent
   Index Method}, Int. J. of Nonlinear Sci., 7, 3, 312-318,
   1579-1585, (2009).

   \item[[14]]. Xue-dong Zheng, Tie-cheng Xia and Hong-qing Zhang, "New Exact Traveling Wave Solutions for Compound
   KdV-Burgers Equations in Mathematical Physics", Appl. Math. E-Notes,
   2, 45-50, (2002).

   \item[[15]]. A.R. Osborne, E. Segre "Numerical solutions of the Korteweg-de Vries equation
   using the periodic scattering transform $\mu$-representation", 44, 3, 575–604, (1990).

   \item[[16]]. Ahmad Farouk, Ekhlass S. Al-Rawi, "Numerical Solution for Non-linear Korteweg-de
   Vries-Burger's Equation Using the Haar Wavelet Method", Iraqi Journal of Statistical Science, 20, (2011).

   \item[[17]]. F. Stenger, "Numerical methods based on sinc and analytical functions", Springer, New York, 1993.

   \item[[18]]. Siraj-ul-Islam, Sirajul Haq, Arshed Ali, "Meshfree method for the numerical solution of the
   RLW equation", Journal of Computational and Applied Mathematics,
   223, 997–1012, (2009).

\end{description}

\end{document}